# AN ENDOMORPHISM OF A FINITELY GENERATED RESIDUALLY FINITE GROUP

DANIEL T. WISE

ABSTRACT. Let $\phi : G \to G$ be an endomorphism of a finitely generated residually finite group. R. Hirshon asked if there exists $n$ such that the restriction of $\phi$ to $\phi^n(G)$ is injective. We give an example to show that this is not always the case.

## 1. INTRODUCTION

A well known theorem of Malćev [LS] states that any surjective endomorphism of a finitely generated residually finite group is an automorphism. In [Hi], R. Hirshon gave the following interesting generalization of Malćev's theorem:

**Theorem 1** (Hirshon). *Let $G$ be a finitely generated residually finite group, And let $\phi$ be an endomorphism onto a finite index subgroup. Then there exists $n$ such that $\phi$ restricted to $\phi^n(G)$ is an injection. In particular, if $G$ is torsion-free, then $\phi$ is an injection.*

A special case of theorem 1 was proven later but independently in [Si] for knot groups.

Hirshon posed the following problem in [Hi] concerning a tempting generalization of theorem 1. We found this problem in the collection of open problems in [M].

**Problem 1.** *Let $G$ be a finitely generated residually finite group and let $\phi$ be an endomorphism of $G$. Does there exist $n$, such that the restriction of $\phi$ to $\phi^n(G)$ is injective?*

The purpose of this paper is to give an example illustrating that the answer to Hirshon's question is no. But before proceeding to give the example in section 2, we wish to pose a problem of our own.

In [Se], Z. Sela proved the impressive result that torsion-free, word-hyperbolic groups are Hopfian. We therefore pose the following variation on problem 1:

**Problem 2.** *Let $G$ be word-hyperbolic, and let $\phi$ be an endomorphism of $G$. Does there exist $n$, such that the restriction of $\phi$ to $\phi^n(G)$ is injective?*

## 2. THE EXAMPLE

It seems clearer to explain the example within a more general context.

We will use the following lemma to show that the endomorphism of our example has the desired property.

*Date*: September 22, 1997.

1991 *Mathematics Subject Classification.* 20E26.

Supported by NSF grant no. DMS-9627506.

I am grateful to Jon McCammond for helpful comments.

**Lemma 1.** *Let $F *_{H=H'} F'$ denote the double of $F$ along the subgroup $H$. Let $\phi : F \to F$ be an endomorphism of $F$ such that $\phi(H) \subset H$. Because $\phi(H) \subset H$, there is an obvious induced endomorphism $\Phi : F *_{H=H'} F' \to F *_{H=H'} F'$. Suppose that $\phi^{-n}(H)$ is a stricly increasing sequence of subgroups of $F$, then $Kernel(\Phi^n)$ is a strictly increasing sequence of subgroups.*

*Proof.* By hypothesis, for each $n$ there exists an element $g_n$ such that $\phi^{n-1}(g_n) \notin H$ but $\phi^n(g_n) \in H$. We conclude that $Kernel(\Phi^n)$ is a strictly increasing sequence of subgroups, because by the normal form theorem, $g_n^{-1} g'_n \in Kernel(\Phi^n)$ but $g_n^{-1} g'_n \notin Kernel(\Phi^{n-1})$. ☐

We will use the following lemma to show that the group of our example is residually finite. Recall that a subgroup $H \subset F$ is said to be *closed* (in the profinite topology) provided that $H$ is the intersection of a set of finite index subgroups of $F$.

**Lemma 2.** *Let $H$ be a closed subgroup of the residually finite group $F$, then the double $F *_{H=H'} F'$ is residually finite.*

*Proof.* By the normal form theorem, an element $w$ of $F *_{H=H'} F'$ is non-trivial if and only if either $w \in H - 1_F$ or $w$ can be represented as a product $w = g_1 g_2 \ldots g_n$ of elements alternately in $F - H$ and $F' - H'$.

In the first case, we consider the obvious homomorphism $F *_{H=H'} F \to F$ and note that $w$ is mapped to a non-trivial element in $F$, and consequently $w$ can be mapped to a non-trivial element in some finite quotient of $F$.

In the second case, we use the hypothesis that $H$ is a closed subgroup to choose a finite index subgroup $K \subset F$ such that for each $i$, $g_i \notin K$. We let $F \to \bar{F}$ be the finite quotient obtained by letting $F$ act on the left cosets of $K$. Let $\bar{H}$ denote the images of $H$. and consider the induced map $F *_{H=H'} F' \to \bar{F} *_{\bar{H}=\bar{H}'} \bar{F}'$.

Let $\bar{g}_i$ denote the image of $g_i$ in $\bar{F}$ and oberserve that $\bar{g}_i \notin \bar{H}$. It follows that the image of $w$, which is the product $\bar{g}_1 \bar{g}_2 \ldots \bar{g}_n$ is a non-trivial normal form in $\bar{F}_{\bar{H}=\bar{H}'} \bar{F}'$. Since $\bar{F}_{\bar{H}=\bar{H}'} \bar{F}'$ is an amalgamated free product of finite groups, it is well known to be residually finite [B]. It follows that $w$ survives in some finite quotient as claimed. ☐

FIGURE 1. The core of $\hat{B}$

**Theorem 2.** *There exists an endomorphism $\Phi : G \to G$ of a finitely generated, residually finite group such that for each $n$, the restriction of $\Phi$ to $\Phi^n(G)$ is not injective.*

*Proof.* Combining lemma 1 and lemma 2 we see that it is sufficient to give an example of (1) a finitely generated, residually finite group $F$, (2) an endomorphism $\phi : F \to F$, and (3) a closed $\phi$-invariant subgroup $H$ such that $\phi^{-n}(H)$ is strictly increasing.

Combining lemma 1 and lemma 2 we see that it is sufficient to give an example of

1. A finitely generated, residually finite group $F$
2. An endomorphism $\phi : F \to F$.
3. A closed $\phi$-invariant subgroup $H$ such that $\phi^{-n}(H)$ is strictly increasing.

We then let $G = F *_{H=H'} F'$ be the double of $F$ along $H$, and let $\Phi : G \to G$ be the endomorphism induced by $\phi : F \to F$ and $\phi : F' \to F'$.

One such example goes as follows: Let $F = \langle a, t \rangle$ be a rank 2 free group. Let $\phi : F \to F$ be induced by $a \mapsto a^2$ and $t \mapsto t$. Let $H = \langle (a^{2^m})^{t^n} : m \geq 0 \rangle$. (We use the notation $x^y = yxy^{-1}$.) It is clear that $F$ is finitely generated and residually finite [Ha]. To see that $\phi(H) \subset H$, observe that $\phi$ maps each generator of $H$ to its square.

To see that $\phi^{-n}(H)$ is strictly increasing we argue as follows: Let $g_n = a^{t^n}$ and observe that $\phi^n(g_n) = \phi^n(a^{t^n}) = (a^{2^n})^{t^n} \in H$, but $\phi^{n-1}(g_n) = \phi^{n-1}(a^{t^n}) = (a^{2^{n-1}})^{t^n} \notin H$. An easy way to see that $(a^{2^{n-1}})^{t^n} \notin H$ is to consider the based covering space $\hat{B} \to B$ (described below) which corresponds to $H$ and to observe that $(a^{2^{n-1}})^{t^n}$ does not lift to a closed path at the basepoint of $\hat{B}$.

To see that $H$ is a closed subgroup, we first note that the finitely generated subgroups of a free group are closed [Ha], and of course the intersection of any set of closed subgroups is closed. It is therefore sufficient to show that $H$ is the intersection of finitely generated subgroups of $F$. Let $H_r = \langle (a^{2^m})t^m, t^{r+1} : 0 \leq m \leq r \rangle$. We will show that $H = \cap_r H_r$.

For each $r$, we have $H \subset H_r$. To see this we show that each generator of $H$ is an element of $H_r$. Consider the generator $t^m a^{2^m} t^{-m}$, and express $m = a(r+1) + b$ where $0 \leq b < (r+1)$ and $a \geq 0$. Then $t^m a^{2^m} t^{-m} = t^{a(r+1)+b} a^{2^{a(r+1)+b}} t^{-(a(r+1)+b)} = (t^{(r+1)})^a (t^b a^{2^b} t^{-b})^{2^{a(r+1)}} (t^{(r+1)})^{-a}$ where the elements in parenthesis are generators of $H_r$.

On the other hand, let $g \notin H$, then we show that there exists $r$, such that $g \notin H_r$.

It seems easiest to prove this using covering spaces. Let $B$ denote a bouquet of two oriented circles labeled $a$ and $t$. Let $\hat{B}$ denote the based covering space of $B$ corresponding to $H$. First consider a based copy of the universal cover of the $t$ circle. Then for each $j \geq 0$, add a copy of the degree $2^j$ cover of the $a$ circle added at the end of the path $t^j$ beginning at the basepoint. Finally, complete this labelled graph to a covering space by adding trees to the incomplete vertices in the obvious way. The core of this based covering space is partially illustrated in figure 2. The basepoint is indicated by a large vertex, and the edges labelled by black and white arrows correspond to the letters $a$ and $t$ respectively. Similarly, we form the based cover $\hat{B}_r$ corresponding to $H_r$, by beginning with the degree $r + 1$ cover of the $t$ circle, and for $0 \leq j \leq r$ we add a degree $2^j$ cover of the $a$ circle at the end of the path $t^j$ beginning at the basepoint. Finally we complete this labelled graph to a covering space of $B$ by adding trees to the incomplete vertices in the obvious way. The core of the based covering space $\hat{B}_3$ is illustrated in figure 2. Note that we follow the same conventions as indicated for figure 2.

Now observe that the ball of radius $n$ about the basepoint of $\hat{B}$ and $\hat{B}_{n+1}$ are isomorphic as labelled graphs. Now suppose that $g \notin H$, let $w$ denote a reduced word representing $g$, then $w$ does not lift to a closed path in $\hat{B}$. Let $|w|$ denote the length of $w$, then since $\hat{B}$ and $\hat{B}_{|w|+1}$ have isomorphic ball's of radius $|w|$, we see that $w$ does not lift to a closed path in $\hat{B}_{|w|+1}$ and hence $g \notin B_{|w|+1}$ and we are done.                                                                                                              $\square$

FIGURE 2. The core of $\hat{B}_3$

DEPARTMENT OF MATHEMATICS, CORNELL UNIVERSITY, ITHACA, NY 14853
*E-mail address*: `daniwise@math.cornell.edu`